\documentclass{amsart}
\usepackage{cases}
\usepackage{amsmath}
\usepackage{amssymb}
\usepackage{amsfonts}
\newtheorem{theorem}{Theorem}[section]
\newtheorem{lemma}[theorem]{Lemma}
\newtheorem{corollary}[theorem]{Corollary}
\newtheorem{proposition}[theorem]{Proposition}

\theoremstyle{definition}
\newtheorem{definition}[theorem]{Definition}
\newtheorem{example}[theorem]{Example}

\newtheorem{remark}[theorem]{Remark}
\numberwithin{equation}{section}
\begin{document}
\title[A note on the splitting theorem for the weighted measure]% end with percent
{A note on the splitting theorem for the weighted measure}
\author{Jia-Yong Wu}

\address{Department of Mathematics, Shanghai Maritime University,
Haigang Avenue 1550, Shanghai 201306, P. R. China}

\email{jywu81@yahoo.com}

\thanks{This work is partially supported by the NSFC11101267.}


\date{\today}

\dedicatory{}

\begin{abstract}
In this paper we study complete manifolds equipped with smooth measures
whose spectrum of the weighted Laplacian has an optimal positive lower
bound and the $m$-dimensional Bakry-\'{E}mery Ricci curvature is
bounded from below by some negative constant. In particular, we prove
a splitting type theorem for complete smooth measure manifolds that
have a finite weighted volume end. This result is regarded as a
study of the equality case of an author's theorem (J. Math. Anal. Appl.
361 (2010) 10-18).
\end{abstract}
\maketitle

\section{Introduction and main result}
The splitting phenomenon for complete manifolds is an interesting topic in
geometric analysis. Perhaps the most notable result is the work of
Cheeger and Gromoll \cite{[Che-Gro1],[Che-Gro2]}, where they proved that
if an $n$-dimensional complete manifold $M$ with nonnegative Ricci
curvature has a geodesic line, then it is isometric to
$\mathbb{R}\times N$ with the product metric, where $N$ is
an $(n-1)$-dimensional complete manifold with nonnegative Ricci curvature.
In a recent work of Wang \cite{[Wang2]}, he proved a splitting theorem
for complete smooth measure manifolds whose $m$-dimensional
Bakry-\'{E}mery Ricci curvature is bounded from below by a negative
multiple of the lower bound of the weighted spectrum. In particular,
from Wang's result, we have
\begin{theorem}\label{Rmain}
Let $(M,g)$ be an $n$-dimensional ($n\geq3$) complete Riemannian manifold
and $\varphi$ be a smooth function on $M$. Assume that the $m$-dimensional
($m\geq n$) Bakry-\'{E}mery Ricci curvature satisfies
\[
Ric_{m,n}\geq-(m-1).
\]
Let $\lambda_1(M)$ be the lower bound of the spectrum of the weighted
Laplacian $\Delta_\varphi=\Delta-\nabla\varphi\cdot\nabla$ on $M$, and
assume that
\[
\lambda_1(M)\geq(m-2).
\]
Then either

\begin{enumerate}
  \item $M$ has only one end with infinite weighted volume; or

\item $M=\mathbb{R}\times N$ with the warped product metric
\[
ds_M^2=dt^2+\cosh^2tds_N^2,
\]
where $N$ is an $(n-1)$-dimensional compact Riemannian manifold.
In this case, $\lambda_1(M)=m-2$.
\end{enumerate}
\end{theorem}
Theorem \ref{Rmain} generalized the work of Li-Wang \cite{[Li-Wang1]} on
Riemannian manifolds to the weighted measure case. If $\varphi$ is constant,
then $Ric_{m,n}=Ric$ for all $m(\geq n)$ and Theorem \ref{Rmain}
returns to the Li-Wang's theorem \cite{[Li-Wang1]} by taking $m=n$.

The weighted measure concept, used in Theorem \ref{Rmain}, can be
briefly described as follows. Let $(M,g)$ be an $n$-dimensional
complete Riemannian manifold and $\varphi$ be a smooth function.
We may define the weighted Laplacian
\[
\Delta_\varphi:=\Delta-\nabla \varphi\cdot\nabla,
\]
which is the infinitesimal generator of the Dirichlet form
\[
\mathcal {E}(\phi_1,\phi_2)=\int_M \langle\nabla \phi_1, \nabla \phi_2\rangle d\mu,
\,\,\,\forall \phi_1, \phi_2\in C_0^{\infty}(M),
\]
where $\mu$ is an invariant measure of $\Delta_\varphi$ given by $
d\mu=e^{-\varphi}dv(g)$. The weighted Laplacian $\Delta_\varphi$
is self-adjoint with respect to the weighted measure $d\mu$.
For the smooth metric measure manifold $(M,g,e^{-\varphi}dv)$, we
define the $m$-dimensional Bakry-\'{E}mery Ricci curvature
(see \cite{[BE],[BQ1],[BQ2],[LD]}) by
\[
Ric_{m,n}:=Ric +Hess(\varphi)-\frac{\nabla \varphi \otimes \nabla
\varphi}{m-n},
\]
where $Ric$ and $Hess$ denote the Ricci curvature and the Hessian of
the metric $g$, respectively. Here $m:=\mathrm{dim}_{BE}(\Delta_\varphi)\geq
n$ is called the Bakry-\'{E}mery dimension of $\Delta_\varphi$, which is a
constant, and $m=n$ if and only if $\varphi$ is a constant \cite{[LD],[Lott]}.
A remarkable feather of $Ric_{m,n}$ is that the volume comparison theorems
hold for $Ric_{m,n}$ in $(M^n,g,e^{-\varphi}dv)$ that look like the case
of Ricci tensor in a $m$-dimensional complete manifold \cite{[LD],WeiWy}.

If we let $m$ be infinite, then the $m$-dimensional Bakry-\'{E}mery Ricci
curvature becomes the $\infty$-dimensional Bakry-\'{E}mery Ricci
curvature
\[
Ric_\infty:=\lim_{m\to \infty}Ric_{m,n}=Ric+Hess(\varphi).
\]
This curvature is closely related to the gradient Ricci soliton
\[
Ric_\infty=\rho g
\]
for some constant $\rho$, which plays an important role in the
theory of Ricci flow \cite{[Cao1]}.

Recently, Fang, Li and Zhang \cite{[FLZ]} obtained two generalizations
of Cheeger-Gromoll splitting theorem via the Bakry-\'{E}mery Ricci
curvature. Munteanu and Wang \cite{[MuWa]} studied function theoretic
and spectral properties on complete noncompact smooth metric measure space
with the nonnegative $\infty$-dimensional Bakry-\'{E}mery Ricci
curvature. In particular, they obtained an interesting splitting
result on complete noncompact gradient steady Ricci solitons.

Using the classical trick of deriving gradient estimates, which is
originated by Yau \cite{[Yau]} (see also \cite{[Cheng-Yau],[SchYau]})
the author proved the following result by choosing
$K=\frac{m-1}{n-1}$ in Theorem 2.1 of \cite{[Wu]}.
\begin{theorem}\label{wthm}
Let $(M,g)$ be an $n$-dimensional ($n\geq2$) complete Riemannian
manifold and $\varphi$ be a smooth function on $M$. Assume that the
$m$-dimensional Bakry-\'{E}mery Ricci curvature satisfies
\[
Ric_{m,n}\geq-(m-1).
\]
Then
\[
\lambda_1(M)\leq\frac{(m-1)^2}{4}.
\]
Moreover, if $f$ be a positive function satisfying
\[
\Delta_\varphi f=-\lambda f
\]
for some constant $\lambda\geq 0$, then $f$ must satisfy the
gradient estimate
\[
|\nabla \ln f|^2\leq\frac{(m-1)^2}{2}-\lambda
+\sqrt{\frac{(m-1)^4}{4}-(m-1)^2\lambda}.
\]
\end{theorem}
Theorem \ref{wthm} was also independently proved by Wang
\cite{[Wang1]}. This result can be viewed as a weighted
version of Cheng's theorem \cite{[Cheng]}. For the case of
$Ric_\infty$, if $|\nabla \varphi|$ is bounded, then we have
another version of gradient estimates \cite{[Wu1]}. Moreover,
the following example shows that the above gradient estimate 
is sharp.

\begin{example}
Consider the $n$-dimensional complete manifold
$M=\mathbb{R}\times N$ endowed with the warped product metric
\[
ds_M^2=dt^2+\exp(-2t)ds_N^2.
\]
If $\{\bar{e}_\alpha\}$ for $\alpha=2,...,n$ form an orthogonal basis of
the tangent space of $N$, then $e_1=\frac{\partial}{\partial t}$ together
with $\{e_\alpha=\exp(-t)\bar{e}_\alpha\}$ form an orthogonal basis for
the tangent space of $M$. By the routine computation, we have
\[
Ric_{M,1j}=-(n-1)\delta_{1j}
\]
and
\[
Ric_{M,\alpha\beta}=\exp(2t)Ric_{N,\alpha\beta}-(n-1)\delta_{\alpha\beta}.
\]
If we choose the weighted function
\[
\varphi:=(m-n)t,
\]
then the $m$-dimensional Bakry-\'{E}mery Ricci curvature of $M$ is
\begin{equation*}
\begin{aligned}
Ric_{mn,1j}&=Ric_{M,1j}+\varphi_{1j}-\frac{\varphi_1\varphi_j}{m-n}\\
&=-(m-1)\delta_{1j}
\end{aligned}
\end{equation*}
and
\[
Ric_{mn,\alpha\beta}=\exp(2t)Ric_{N,\alpha\beta}-(n-1)\delta_{\alpha\beta}.
\]
Hence we observe that if the Ricci curvature of manifold $N$
is nonnegative, then the $m$-dimensional Bakry-\'{E}mery Ricci
curvature satisfies
\[
Ric_{mn}\geq-(m-1).
\]
In this setting, we \emph{claim} that
\[
\lambda_1(M)=\frac{(m-1)^2}{4}.
\]
Indeed, we may choose the function
\[
f=\exp\left(\frac{m-1}{2}t\right).
\]
A direct computation yields that
\begin{equation*}
\begin{aligned}
\Delta_\varphi f&=\frac{d^2f}{dt^2}-(n-1)\frac{df}{dt}
-\frac{d\varphi}{dt}\cdot\frac{df}{dt}\\
&=-\frac{(m-1)^2}{4}f,
\end{aligned}
\end{equation*}
since $\Delta=\frac{\partial^2}{\partial t^2}
-(n-1)\frac{\partial}{\partial t}+\exp(2t)\Delta_N$.
On the other hand, we have the following proposition, which is a mild
generalization for the classical case.
\begin{proposition}\label{pro}
Let $(M,g)$ be an $n$-dimensional complete Riemannian
manifold and $\varphi$ be a smooth function on $M$.
If there exists a positive function $f$ satisfying
\[
\Delta_\varphi f\leq-\lambda f,
\]
then $\lambda_1(M)$, the lower bound of the spectrum of the weighted
Laplacian $\Delta_\varphi$, satisfies
\[
\lambda_1(M)\geq \lambda.
\]
\end{proposition}
Combining this proposition and Theorem \ref{wthm}, we immediately
conclude that $\lambda_1(M)=\frac{(m-1)^2}{4}$ as claimed.
\end{example}

\vspace{0.5em}

Since $\frac{(m-1)^2}{4}\geq m-2$ with equality holds only when $m=3$
(in this case, we return to the classical Laplacian case,
see Remark \ref{rem}), Theorem \ref{Rmain} in fact asserts that
the equality case in Theorem \ref{wthm} implies that the measure
manifold belongs to the case (1) of Theorem \ref{Rmain}.
Namely, the measure manifold must only have one infinite
weighted volume end, unless $m=3$. In this case, the warped
product given in Theorem  \ref{Rmain} is the only exception.

Naturally, we would like to ask if the finite weighted volume ends can
be ruled out when a measure manifold satisfies the hypotheses of
Theorem \ref{Rmain}. In this paper, we follow the arguments of 
Li-Wang's work \cite{[Li-Wang2]}, and show that the above example 
is the only case (it may be different from the weighted function 
$\varphi$) when $M$ has a finite weighted volume end if $M$ achieves
equality in weighted spectrum upper bound of Theorem \ref{wthm}.
\begin{theorem}\label{main}
Let $(M,g)$ be an $n$-dimensional ($n\geq3$) complete Riemannian
manifold and $\varphi$ be a smooth function on $M$.
Assume that the $m$-dimensional ($m>3$)
Bakry-\'{E}mery Ricci curvature satisfies
\[
Ric_{m,n}\geq-(m-1).
\]
and $\lambda_1(M)$, the lower bound of the spectrum of the weighted
Laplacian $\Delta_\varphi$, satisfies
\[
\lambda_1(M)\geq\frac{(m-1)^2}{4}.
\]
Then either
\begin{enumerate}
  \item $M$ has only one end; or

\item $M=\mathbb{R}\times N$ with the warped product metric
\[
ds_M^2=dt^2+\exp(-2t)ds_N^2,
\]
where $N$ is an $(n-1)$-dimensional compact manifold. Moreover,
\[
\varphi(t,x)=\varphi(0,x)+(m-n)t
\]
for all $(t,x)\in\mathbb{R}\times N$.
\end{enumerate}
\end{theorem}
\begin{remark}\label{rem}
In Theorem \ref{main}, we assume that $m>3$.  Because when $m=3$, we observe that
$\frac{(m-1)^2}{4}=m-2$, $n=m=3$ and hence $\varphi$ is constant.
Therefore, this case is exact the Li-Wang classical
result (Theorem 0.6 in \cite{[Li-Wang2]}).
\end{remark}

\begin{remark}
Using similar trick, we can obtain splitting theorems for complete
manifolds with $\infty$-dimensional Bakry-\'{E}mery curvature
by gradient estimates of \cite{[Wu1]}.
This was treated by the author in a separated paper \cite{[Wu2]}.
\end{remark}

% ------------------------------------------------------------------------

\section{Preliminary}
In this section, we will give some important lemmas to prepare the
proof of Theorem \ref{main}. At first, we recall some basic
definitions in smooth metric measure manifolds, which are also introduced
in \cite{[Wang2]}.
\begin{definition}
Let $(M,g)$ be a complete Riemannian manifold and $\varphi$ be a smooth
function on $M$. A weighted Green's function $G_\varphi(x,y)$ is a function
defined on $(M\times M)\backslash\{(x,x)\}$ satisfying
\begin{enumerate}
  \item $G_\varphi(x,y)=G_\varphi(y,x)$, and

\item $\Delta_{\varphi,y}G(x,y)=-\delta_{\varphi,x}(y)$,
\end{enumerate}
for all $x\neq y$, where $\delta_{\varphi,x}(y)$ is defined by
\[
\int_M\psi(y)\delta_{\varphi,x}(y)d\mu=\psi(x)
\]
for every compactly supported function $\psi$.
\end{definition}
In fact, every smooth measure manifold admits a weighted Green's function.
Following Li-Tam \cite{[Li-Tam]}, we can give a constructive argument for the
existence of $G_\varphi(x,y)$. But some measure manifolds admit weighted Green's
functions which are positive and others may not. This special property
distinguishes the weighted function theory of complete measure
manifolds into two classes.
\begin{definition}
A complete measure manifold $(M,g,e^{-\varphi}dv)$ is said to be weighted
non-parabolic if it admits a positive weighted Green's function.
Otherwise, it is said to be weighted parabolic.
\end{definition}

Following the arguments of Theorem 2.3 in \cite{[PLi]}, we can easily
show that a complete measure manifold is weighted non-parabolic if and
only if there exists a positive weighted super-harmonic function whose
infimum is achieved at infinitely. In the following, we will give the
definition of an end of a complete manifold.
\begin{definition}
An end, $E$, with respect to a compact subset
 $\Omega\subset M$ is an unbounded connected component
of $M\setminus \Omega$. The number of ends with respect of
$\Omega$, denoted by $N_\Omega(M)$, is the number of unbounded connected
component of $M\setminus \Omega$.
\end{definition}
It is easy to see that if $\Omega_1\subset \Omega_2$, then $N_{\Omega_1}(M)
\leq N_{\Omega_2}(M)$. Hence if $\Omega_i$ is a compact exhaustion of $M$,
then $N_{\Omega_i}(M)$ is a monotonically nondecreasing sequence. If this
sequence is bounded, then we say that $M$ has finitely many ends. In this case,
the number of ends of $M$ is defined by
\[
N(M)=\lim_{i\to\infty} N_{\Omega_i}(M).
\]
Obviously, the number of ends is independent of the compact exhaustion
$\{\Omega_i\}$.

\begin{definition}
An end $E$ is said to be weighted non-parabolic if it admits a
positive weighted Green's function with Neumann boundary
condition on $\partial E$. Otherwise, it is said to be
weighted parabolic.
\end{definition}
From the construction of  Li-Tam \cite{[Li-Tam]}, we can easily
verify that a complete measure manifold is weighted
non-parabolic if and only if it has a weighted non-parabolic end.

We now state a decay property about weighted
harmonic functions on the end of a smooth metric measure manifold,
which is a slight generalization of Lemma 1.1 in \cite{[Li-Wang1]}.
\begin{lemma}\label{lem0}
Let $(M,g)$ be an $n$-dimensional complete Riemannian manifold and
$\varphi$ be a smooth function on $M$. Suppose $E$ is an end of $M$ and the
weighted spectrum $\lambda_1(M)>0$. Then for any weighted harmonic
function $f$ on $E$ such that $f=\lim_{R_i\to\infty}f_i$ with
$\Delta_\varphi f_i=0$ on $E(R_i)$ and $f_i=0$ on
$E\cap \partial B_p(R_i)$, $f$ satisfies the decay estimate
\[
\int_{E(R+1)\setminus E(R)}f^2d\mu\leq C\exp(-2\sqrt{\lambda_1(M)}R)
\]
for some constant $C>0$ depending on $f$, $\lambda_1(M)$ and $n$,
where $B_p(R)$ denotes a geodesic ball centered at some fixed
point $p\in M$ with radius $R>0$, and $E(R)=B_p(R)\cap E$.
\end{lemma}

The following lemma is an characterization
for an end by its weighted volume due to Wang \cite{[Wang2]}.
\begin{lemma}\label{lem1}
Let $(M,g)$ be an $n$-dimensional complete Riemannian manifold and
$\varphi$ be a smooth function on $M$.  We assume that
\[
\lambda_1(M)\geq\frac{(m-1)^2}{4}.
\]
Let $E$ be an end of $M$, and let $V_\varphi(E)$ be the simply
weighted volume of end $E$. We denote the weighted volume of the
set $E(R)$ by $V_\varphi(E(R))$. $R>0$ is large enough.
\begin{enumerate}
  \item  If $E$ is a weighted parabolic end, then E must have exponential
weighted volume decay given by
\[
V_\varphi(E)-V_\varphi(E(R))\leq C\exp(-(m-1)R)
\]
for some constant $C>0$ depending on the end $E$.

\item If $E$ is a weighted non-parabolic end, then $E$ must
have exponential volume growth given by
\[
V_\varphi(E(R))\geq C\exp((m-1)R)
\]
for some constant $C>0$ depending on the end $E$.
\end{enumerate}
\end{lemma}
\begin{remark}
Lemma \ref{lem1} can be viewed as a refined version of
Theorem \ref{wthm}. In fact, if the
$m$-dimensional Bakry-\'{E}mery Ricci curvature satisfying
\[
Ric_{m,n}\geq-(m-1),
\]
then the weighted Bishop volume comparison theorem
(see \cite{[BQ2],[LD]}) asserts that
\[
V_\varphi(B_p(R))\leq V_{\mathbb{H}^m}(B_p(R))\leq C\exp((m-1)R).
\]
Combining this  and Lemma \ref{lem1}, we conclude that
\[
\lambda_1(M)\leq\frac{(m-1)^2}{4},
\]
as asserted in Theorem \ref{wthm}.
\end{remark}

On the other hand, if the $m$-dimensional Bakry-\'{E}mery Ricci
curvature is bounded from below by $-(m-1)$, then the weighted
Bishop volume comparison theorem, says that for any $x\in M$,
\[
\frac{V_\varphi(B_x(R))}{V_{\mathbb{H}^m}(B(R))}
\]
is nonincreasing in $R$, where
$V_\varphi(B_x(R))=\int_{B_x(R)}e^{-\varphi}dv(g)$ denotes the weighted
volume of the geodesic ball $B_x(R)$, and $V_{\mathbb{H}^m}(B(R))$
denotes the volume of a geodesic ball of radius $R$ in the $m$-dimensional
hyperbolic space form $\mathbb{H}^m$ with constant curvature $-1$.
Therefore for any $R_1<R_2$, we have
\[
\frac{V_\varphi(B_x(R_2))}{V_\varphi(B_x(R_1))}
\leq\frac{V_{\mathbb{H}^m}(B(R_2))}{V_{\mathbb{H}^m}(B(R_1))}.
\]
In particular, if we let $x=p$, $R_1=0$ and $R_2=R$, then
\begin{equation}\label{bijiao1}
V_\varphi(B_p(R))\leq C\exp((m-1)R)
\end{equation}
for sufficiently large $R$. If we let $x\in\partial B_p(R)$, $R_1=1$
and $R_2=R+1$, then
\begin{equation}
\begin{aligned}\label{bijiao2}
V_\varphi(B_x(1))&\geq CV_\varphi(B_x(R+1))\exp(-(m-1)R)\\
&\geq CV_\varphi(B_p(1))\exp(-(m-1)R).
\end{aligned}
\end{equation}
Combining \eqref{bijiao1}, \eqref{bijiao2} and Lemma \ref{lem1},
we have that
\begin{corollary}\label{lem2}
Let $(M,g)$ be a complete Riemannian manifold and $\varphi$ be a smooth
function on $M$, with the $m$-dimensional Bakry-\'{E}mery Ricci
curvature satisfying
\[
Ric_{m,n}\geq-(m-1).
\]
We assume that
\[
\lambda_1(M)\geq\frac{(m-1)^2}{4}.
\]
Let $E$ be an end of $M$, and let $V_\varphi(E)$ be the simply
weighted volume of end $E$. We denote the weighted volume of the
set $E(R)$ by $V_\varphi(E(R))$. $R>0$ is large enough.
\begin{enumerate}
  \item  If $E$ is a weighted-parabolic end, then E must have exponential
weighted volume decay given by
\[
C_4\exp(-(m-1)R)\leq V_\varphi(E)-V_\varphi(E(R))\leq C_1\exp(-(m-1)R)
\]
for some constant $C_1\geq C_4>0$ depending on the end $E$.

\item If $E$ is a weighted-non-parabolic end, then $E$ must
have exponential volume growth given by
\[
 C_3\exp((m-1)R)\geq V_\varphi(E)\geq C_2\exp((m-1)R)
\]
for some constant $C_3\geq C_2>0$ depending on the end $E$.
\end{enumerate}
\end{corollary}

% ------------------------------------------------------------------------

\section{Proof of Theorem \ref{main}}
We are now ready to prove Theorem \ref{main} in introduction. The
proof method belongs to Li-Wang \cite{[Li-Wang2]}.
\begin{proof}[Proof of Theorem \ref{main}]
Suppose that the manifold $M$ satisfies the hypothesis of
Theorem \ref{main}. Then Theorem \ref{Rmain} asserted that $M$ must have only
one infinite weighted volume end because the warped product with the metric
given by
\[
ds_M^2=dt^2+\cosh^2tds_N^2
\]
has $\lambda_1(M)=m-2$, which does not satisfy the second hypothesis of
Theorem \ref{main}.

Now we assume that manifold $M$ has a finite weighted volume end.
Since $\lambda_1(M)>0$, $M$ must also have an infinite weighted
volume end. By choosing the compact set $D$ appropriately,
we may assume that $M\backslash D$ has one infinite weighted volume,
weighted non-parabolic end $E_1$ and one finite weighted volume,
weighted parabolic end $E_2$.

In an analogous way as Li-Tam's arguments \cite{[Li-Tam]},
our consideration is the weighted measure case. We assert that
there exists a positive weighted harmonic function
$f$ with the the following properties:

\begin{itemize}
\item $\inf_{\partial E_1(R)}\to 0$   as $R\to\infty$;

\item $\sup_{\partial E_2(R)}\to \infty$   as $R\to\infty$; and

\item $f$ is bounded and has finite weighted Dirichlet integral on $E_1$.
\end{itemize}
Then the gradient estimate of Theorem \ref{wthm} implies that
\[
|\nabla f|^2\leq (m-1)^2f^2.
\]
Combining this with the fact that function $f$ is weighted harmonic, we have
\begin{equation}
\begin{aligned}\label{jisuan0}
\Delta_\varphi f^{1/2}&=-\frac 14f^{-3/2}|\nabla f|^2\\
&\leq-\frac{(m-1)^2}{4}f^{1/2}.
\end{aligned}
\end{equation}
If we let $h=f^{1/2}$, then for any nonnegative cut-off function $\psi$ we have
\begin{equation}\label{indent1}
\int_M|\nabla(\psi h)|^2d\mu=\int_M|\nabla\psi|^2 h^2d\mu+\int_M\psi^2|\nabla h|^2d\mu
+2\int_M\psi h \nabla\psi \nabla hd\mu.
\end{equation}
Since
\[
\int_M\psi h \nabla\psi \nabla hd\mu
=-\int_M\psi \nabla\psi h\nabla hd\mu
-\int_M\psi^2 |\nabla h|^2d\mu-\int_M\psi^2h\, \Delta_\varphi hd\mu,
\]
the integral equality \eqref{indent1} reduces to
\begin{equation}
\begin{aligned}\label{jisuan1}
\int_M|\nabla(\psi h)|^2d\mu
&=\int_M|\nabla\psi|^2 h^2d\mu-\int_M\psi^2h\, \Delta_\varphi hd\mu\\
&=\int_M|\nabla\psi|^2 h^2d\mu+\frac{(m-1)^2}{4}\int_M\psi^2h^2\\
&\quad-\int_M\psi^2h\left[\frac{(m-1)^2}{4}h+\Delta_\varphi h\right]d\mu.
\end{aligned}
\end{equation}
Since $\lambda_1(M)\geq\frac{(m-1)^2}{4}$, the definition of $\lambda_1(M)$
gives us
\[
\frac{(m-1)^2}{4}\int_M\psi^2h^2d\mu\leq\int_M|\nabla(\psi h)|^2d\mu.
\]
Hence
\begin{equation}\label{jisuan2}
\int_M\psi^2h\left[\frac{(m-1)^2}{4}h+\Delta_\varphi h\right]d\mu
\leq \int_M|\nabla\psi|^2 h^2d\mu.
\end{equation}

Integrating the gradient estimate of Theorem \ref{wthm} along geodesics,
we know that $f$ must satisfy the growth estimate
\[
f(x)\leq C\exp((m-1)r(x)),
\]
where $r(x)$ is the geodesic distance from $x$ to a fixed point $p\in M$.
In particular, when restricted on the parabolic end $E_2$, together with
the volume estimate of Lemma \ref{lem1}, we conclude that
\begin{equation}\label{jisuan3}
\int_{E_2(R)}fd\mu
\leq CR.
\end{equation}
On the other hand, Lemma \ref{lem0} asserts that on $E_1$, the function
$f$ must satisfy the decay estimate
\[
\int_{E_1(R+1)\setminus E_1(R)}f^2d\mu\leq C\exp(-(m-1)R)
\]
for $R$ sufficiently large. By the Schwarz inequality, we have
\[
\int_{E_1(R+1)\setminus E_1(R)}fd\mu\leq C\exp\left(-\frac{m-1}{2}R\right)
V^{1/2}_{\varphi E_1}(R+1),
\]
where $V_{\varphi E_1}(r)$ denotes the weighted volume of $E_1(r)$.
Combining this with the volume estimate of Corollary \ref{lem2}, we
have that
\[
\int_{E_1(R+1)\setminus E_1(R)}fd\mu\leq C
\]
for some constant $C$ independent of $R$. In particular, we have
\[
\int_{E_1(R)}fd\mu\leq CR.
\]
Combining this with \eqref{jisuan3}, we conclude that
\begin{equation}\label{jisuan4}
\int_{B_p(R)}fd\mu \leq CR.
\end{equation}

Now we define the cut-off function $\psi$
on $M$ in \eqref{jisuan2} by
\begin{equation*}
\psi(x)=\left\{
\begin{aligned}
1 \quad\quad& x\in B_p(R)\\
\frac{2R-r}{R} \quad\quad& x\in B_p(2R)\setminus B_p(R)\\
0 \quad\quad& x\not\in B_p(2R).
\end{aligned}
\right.
\end{equation*}
Hence the right hand side of \eqref{jisuan2} is given by
\[
\int_M|\nabla\psi|^2 h^2d\mu=R^{-2}\int_{B_p(2R)\setminus B_p(R)} h^2d\mu
\]
and \eqref{jisuan4} implies
\[
\int_M|\nabla \psi|^2h^2d\mu \to 0
\]
as $R\to\infty$. Therefore we obtain
\[
\Delta_\varphi h=-\frac{(m-1)^2}{4}h
\]
and inequality \eqref{jisuan0} used in the above argument is an equality. In
particular, we have
\[
|\nabla f|=(m-1)f
\]
and
\begin{equation}\label{equ}
|\nabla (\ln f)|^2=(m-1)^2.
\end{equation}
Hence the inequalities used to prove the gradient estimate of
Theorem \ref{wthm} are all equalities. Namely we must have
equality (2.11) in \cite{[Wu]} since
\[
\Delta_\varphi|\nabla(\ln f)|^2=\Delta|\nabla(\ln f)|^2
-\nabla \varphi\cdot\nabla|\nabla(\ln f)|^2=\nabla|\nabla(\ln f)|^2=0.
\]
Moreover the inequalities used to derive (2.11) in \cite{[Wu]} must all be
equalities. More specifically, equality (2.6) in \cite{[Wu]} implies
\[
(\ln f)_{1j}=0
\]
for all $1\leq j\leq n$, whereas equality (2.7) in \cite{[Wu]} gives
\begin{equation}\label{equality2}
\langle\nabla \varphi,\nabla \ln f\rangle=(m-1)(m-n)
\end{equation}
and
\begin{equation*}
\begin{aligned}
(\ln f)_{\alpha\beta}&=-\frac{|\nabla(\ln f)|^2}{m-1}\delta_{\alpha\beta}\\
&=-(m-1)\delta_{\alpha\beta}
\end{aligned}
\end{equation*}
for all $2\leq \alpha,\beta\leq n$. Since $e_1$ is the unit normal to
the level set of $\ln f$, the second fundamental form $\mathrm{II}$ of
the level set is given by
\begin{equation*}
\begin{aligned}
\mathrm{II}_{\alpha\beta}&=\frac{(\ln f)_{\alpha\beta}}{(\ln f)_1}\\
&=\frac{-(m-1)\delta_{\alpha\beta}}{m-1}\\
&=-\delta_{\alpha\beta}.
\end{aligned}
\end{equation*}
Moreover, \eqref{equ} implies that if we set $t=\frac{\ln f}{m-1}$,
then $t$ must be the distance function between the level sets of $f$,
hence also for $\ln f$. Since
$\mathrm{II}_{\alpha\beta}=(-\delta_{\alpha\beta})$,
this implies that the metric on $M$ can be written as
\[
ds_M^2=dt^2+\exp(-2t)ds_N^2.
\]
By \eqref{equality2}, we also have
\[
\varphi(t,x)=\varphi(0,x)+(m-n)t,
\]
where $(t,x)\in\mathbb{R}\times N$. Since we assume that the
manifold $M$ has two ends, $N$ must be compact.
\end{proof}

% ------------------------------------------------------------------------
\bibliographystyle{amsplain}

\end{document}